\documentclass[english,a4paper]{amsart}
\usepackage{xspace,bm}
\usepackage[dvips,pdftitle=A\ note\ on\ the\ CR\
cohomology\ of\ Levi-Flat\ minimal\ orbits\ in\ complex\ flag\ manifolds,pdfauthor=Andrea\ Altomani]{hyperref}

\title[Cohomology of Levi-flat minimal orbits]
{A note on the CR
cohomology of Levi-Flat minimal orbits in complex flag manifolds}

\author{Andrea Altomani}
\address[A. Altomani]
         {Dipartimento di Matematica\\
         II Universit\`a di Roma ``Tor Vergata''\\
         Via della Ricerca Scientifica\\
         00133 Roma (Italy)}
  \email{altomani@mat.uniroma2.it}

\date{February 28, 2006}

\theoremstyle{plain}
\newtheorem{thm}{Theorem}[section]
\newtheorem{lem}[thm]{Lemma}
\newtheorem{cor}[thm]{Corollary}
\newtheorem{prop}[thm]{Proposition}
\numberwithin{equation}{section}
\newcommand{\CR}{\texorpdfstring{\ensuremath{CR}}{CR}\xspace}
\newcommand{\C}{\mathbb{C}}
\newcommand{\R}{\mathbb{R}}
\renewcommand{\P}{\mathbb{P}}

\newcommand{\fr}{\mathfrak}
\newcommand{\g}{\fr{g}}
\newcommand{\q}{\fr{q}}
\renewcommand{\k}{\fr{k}}
\renewcommand{\l}{\fr{l}}

\newcommand{\gr}{\mathbf}
\newcommand{\G}{\gr{G}}
\newcommand{\Q}{\gr{Q}}
\newcommand{\K}{\gr{K}}
\renewcommand{\L}{\gr{L}}
\newcommand{\N}{\gr{N}}

\newcommand{\sh}{\mathcal}
\renewcommand{\H}{\sh{H}}
\renewcommand{\O}{\sh{O}}
\newcommand{\E}{\sh{E}}

\newcommand{\Z}{\sh{Z}}
\newcommand{\B}{\sh{B}}
\newcommand{\A}{\sh{A}}

\newcommand{\vb}{\bm}
\newcommand{\Eb}{\vb{E}}
\newcommand{\Hb}{\vb{H}}

\newcommand{\e}{\mathrm{e}}

\newcommand{\dbar}{\bar{\partial}}
\newcommand{\dbarNto}{\xrightarrow{\dbar_N}}
\newcommand{\ie}{i.\;e{.}{}\xspace}

\DeclareMathOperator{\Lie}{Lie}
\DeclareMathOperator{\Span}{Span}

\begin{document}
\begin{abstract}
We prove a relation between the $\dbar_M$ cohomology of a minimal orbit
$M$ of a real form $\G_0$ of a complex semisimple Lie group $\G$ in a
flag manifold $\G/\Q$ and the Dolbeault cohomology of the Matsuki dual
open orbit $X$ of the complexification $\K$ of a maximal compact
subgroup $\K_0$ of $\G_0$, under the assumption that $M$ is Levi-flat.
\end{abstract}

\maketitle

\section{Introduction}
Many authors have studied the $\dbar_M$ cohomology of \CR manifolds 
(see e.\;g.\ \cite{HN:2000, HN:2002, LL:2003} and references therein).
In particular, since Andreotti and Fredricks \cite{AF:1979} proved 
that every real analytic \CR manifold $M$ can be embedded in a complex
manifold $X$, it is natural to try to find relations between the
$\dbar_M$ cohomology of $M$ and the Dolbeault cohomology of $X$.
\par
In this paper we examine this problem for a specific class of
homogeneous \CR manifolds, namely  minimal orbits in complex 
flag manifolds that are Levi-flat.
\par
Given a (generalized) flag manifold $Y=\G/\Q$, with $\G$ a complex
semisimple Lie group and $\Q$ a parabolic subgroup of $\G$, a real form
$\G_0$ of $\G$ acts on Y with finitely many orbits.
Among these there is
exactly one orbit that is compact, the minimal orbit $M=\G_0\cdot o$.
Let $\K_0$ be a maximal compact subgroup of $\G_0$, and $\K$ its 
complexification. 
Then $X=K\cdot o$ is an open dense complex submanifold of 
$\G/\Q$ and contains $M$ as an embedded submanifold.
It is known  that $M$ is a deformation retract of $X$,
so $H^p(M,\C)=H^p(X,\C)$ (\cite{CS:1999}, \cite{HW:2003}).
\par
Let $\Eb$ be a $\K$-homogeneous complex vector bundle on $X$ and
$\Eb|_M$
its restriction to $M$.
Under the additional assumption that $M$ is Levi-flat we prove that
the restriction map from the Dolbeault cohomology $H^{p,q}(X,\Eb)$ to
the $\dbar_M$ cohomology $H^{p,q}(M,\Eb|_M)$ is continuous, injective
and
has a dense range.
More precisely we show that
\[
H^{p,q}(M,\Eb|_M)=\O_M(M)\otimes_{\O_X(X)}H^{p,q}(X,\Eb)
\]
where $\O_M(M)$ (resp. $\O_X(X)$) is the space of \CR (resp. 
holomorphic) functions on $M$ (resp. $X$), and that the restriction map
from $\O_X(X)$ to $\O_M(M)$ is injective, continuous and has a dense 
range.
\par

\section{Preliminaries on minimal
orbits 
in complex flag manifolds}
Let $\G$ be a complex connected semisimple Lie group, with Lie algebra
$\g$, and let $\Q$ be a parabolic subgroup of $\G$, with Lie algebra
$\q$.
Then $\Q$ is the normalizer of $\q$, $\Q=\N_{\G}(\q)$ and is connected.
The coset space $Y=\G/\Q$ is a compact complex manifold, called a
\emph{flag manifold} (it is a complex smooth projective variety).
It is not restrictive to assume that $\q$ does not contain any simple
ideal of $\g$.
\par
Let $\sigma$ be an anti-holomorphic involution of $\G$; we will also 
denote by $\sigma$ its differential at the identity
and we will write $\bar x=\sigma(x)$.
A \emph{real form} of $\G$ is an open subgroup $\G_0$ of $\G^\sigma$.
It is a Lie subgroup and its Lie algebra $\g_0$ satisfies 
$\g_0=\g^\sigma$ and $\g=\C\otimes\g_0$.
\par
Let $\K_0$ be a maximal compact subgroup of $\G_0$, and $\theta$ 
the corresponding Cartan involution: $\K_0=\G_0^\theta$.
Still denoting by $\theta$ the complexification of $\theta$, there is
exactly one open subgroup $\K$ of $\G^\theta$ such that $\K\cap\G_0=
\K_0$.
Let $\k$ and $\k_0$ be the corresponding Lie subalgebras.
\par
The groups $\G_0$ and $\K$ act on $Y$ via left multiplication. 
There is exactly one closed $\G_0$-orbit $M$
(\cite{AMN:2005,Wolf:1969}) 
and it 
is contained in the unique open $\K$ orbit $X$
(\cite{BL:2002,Matsuki:1982}), we denote by $j\colon M\to X$ the
inclusion. The open orbit $X$ is dense in $Y$ and is the dual orbit of
$M$, in the terminology of \cite{Matsuki:1982}.
\par
The manifolds $M$  and $X$ do not depend on the choice of $\G$ and 
$\G_0$, but only on $\g_0$ and $\q$ (\cite{AMN:2005, Matsuki:1982}).
So there is no loss of generality assuming that $\G$ is 
simply connected, that
$\G_0$ is connected 
and that $M$ and $X$ are the orbits through the point $o=e\Q$.
We will write $M=M(\g_0,\q)$.
\par
The isotropy subgroups at $o$ of the actions of $\G_0$ and $\K$ are 
$\G_+=\G_0\cap\Q$ and $\L=\K\cap\Q$, with Lie algebras $\g_+=\g_0\cap
\q$ and $\l=\k\cap\q$.
Since $M$ is compact, the action on $M$ of maximal compact subgroup
$\K_0$ is
transitive: $M=\K_0/\K_+$, where $\K_+=\K_0\cap\Q=\G_0\cap\L$ and
$\Lie(\K_+)=\k_+=\k_0\cap\q=\g_0\cap\l$.
\par
In the language of \cite{AMN:2005} the pair $(\g_0,\q)$ is an effective
parabolic minimal \CR algebra and $M$ is the associated minimal orbit.
On $M$ there is a natural \CR structure induced by the inclusion into
$X$.
\par
We recall that $M$ is totally real if the partial complex structure is trivial. 
We give a more complete characterization of totally real minimal orbits:
\begin{thm}\label{thm:totallyreal}
The following are equivalent:
\begin{enumerate}
	\item $M$ is totally real.
	\item $\bar\q=\q$.
	\item $\l=\k_+^\C$.
	\item $X$ is a Stein manifold.
	\item $X$ is a smooth affine algebraic variety defined over $\R$ and $M$ 
	is the set of its real points.
\end{enumerate}
\end{thm} 
\begin{proof}
$(5)\Rightarrow(4)$ 
because every closed complex submanifold of a complex
vector space is Stein.
\par
$(4)\Rightarrow(3)$ 
Since $X$ is Stein, its covering $\tilde
X=\K/\L^0$ is also Stein (\cite{Stein:1956}). 
Furthermore $\K$ is a linear algebraic group that is the 
complexification of a maximal compact subgroup $\K_0$; the
result then follows from Theorem~3 of \cite{Matsushima:1960}
and Remark~2 thereafter.
\par
$(3)\Rightarrow(5)$
If $\l$ is the complexification of $\k_+$, then $\L$ is the
complexification of $\K_+$.
Hence $X=\K/\L$ is the complexification of $M=\K_0/\K_+$ in the sense of
\cite{IS:1966}, and $(5)$ follows from Theorem 3 of the same paper.
\par
$(1)\Leftrightarrow(2)$ 
is easy, and proved in \cite{AMN:2005}.
\par
$(5)\Rightarrow(1)$
is obvious.
\par
$(2)\Rightarrow(3)$
We have that $\bar\k=\k$, thus
$\k\cap\q=(\k\cap\q\cap\g_0)^\C=\k_+^\C$.
\end{proof}

We denote by $\O_N$ the sheaf of smooth \CR functions on a \CR manifold
$N$.
If $N$ is complex or real, then $\O_N$ is the usual sheaf of 
holomorphic or smooth (complex valued) functions.
For every open set $U\subset N$, the space $\O_N(U)$ is a Fr\'echet
space (with the topology of uniform convergence of all derivatives on
compact sets).
\par
\begin{cor}\label{cor:dense} 
If $M$ is totally real then $j^*\big(\O_X(X)\big)$ is dense in $\O_M(M)=C^\infty(M)$.
\end{cor}
\begin{proof}
Let $X\subset\C^N$ be an embedding as in $(5)$ of Theorem~\ref{thm:totallyreal}, so that
$M=X\cap\R^N$.
The restrictions of complex polynomials in $\C^N$ are contained in
$\O_X(X)$ and dense
in $\O_M(M)$ (see e.g{.}\ \cite{Treves:1967}).
\end{proof}

\section{Levi-flat orbits and the fundamental reduction}
In this paper we consider Levi-flat minimal orbits.
They are orbits $M=M(\g_0,\q)$, where
 $\q'=\q+\bar\q$ is a subalgebra 
(necessarily parabolic) of $\g$.
Let $\Q'=\N_{\G}(\q')$, $Y'=\G/\Q'$, $\G_+'=\G\cap\Q'$,
$M'=M(\g_0,\q')=\G_0/\G_+'$, $\K_+'=\K_0\cap\Q'$, $\L'=\K\cap\Q'$ and
$X'=\K/\L'$.
\par
From Theorem \ref{thm:totallyreal} we have that $M'$ is totally real
and $X'$ is Stein.
\par
The inclusion $\Q\to\Q'$ induces a fibration
\begin{equation}
\pi\colon Y=\G/\Q\longrightarrow\G/\Q'=Y'
\end{equation}
with complex fiber $F\simeq\Q'/\Q$.
This fibration is classically referred to as the Levi foliation, and is
a special case of the fundamental reduction of
\cite{AMN:2005}.
In fact Levi-flat minimal orbits are characterized by the property that
the fibers of the fundamental reduction are totally complex.
\par
We identify $F$ with $\pi^{-1}(e\Q)$.

\begin{lem}
$\pi^{-1}(M')=M$, $\pi^{-1}(X')=X$ and $F$ is a compact
connected
complex flag manifold.
\end{lem}
\begin{proof}
First we observe that $F$ is connected because $\Q'$ is connected.
Let $F'$ be the fiber of the restriction $\pi|_M\colon M\to M'$. 
Then $F'$ is totally complex and \CR generic in $F$, hence an open
subset of $F$.
Proposition 7.3 and Theorem 7.4 of \cite{AMN:2005} show
that there exists a connected real semisimple Lie group $\G_0''$
acting on $F'$ with an open
orbit, and that the Lie algebra of the isotropy is a $t$-subalgebra
(\ie contains a maximal triangular subalgebra) of $\g''_0$.
Hence a maximal compact subgroup
$\K_0''$ of $\G_0''$ has an open orbit, which is also closed.
Since $F'$ is open in $F$ and $F$ is connected, $\K_0''$ is transitive
on $F$, and $F'=F$, proving the first two
statements.
\par
Furthermore the isotropy subgroup $\G_+''$ for the action of $\G_0''$
on $F$ and the homogeneous complex structure are exactly those of a
totally complex minimal orbit, hence by \cite[\S\,10]{AMN:2005} $F$ is a
complex flag manifold.
\end{proof}
\par
The total space $M$ is locally isomorphic to an open subset of $M'\times
F$, hence to $U\times\R^k$, where $U$ is open in $\C^n$, for some
integers $n$ and $k$.
\par
For a Levi flat \CR manifold $N$ and a nonnegative integer $p$, let
$\Omega^p_N$ be the sheaf of $p$-forms that are \CR
(see \cite{HN:2002} for precise definitions).
They are $\O_N$-modules and $\Omega^0_N\simeq\O_N$.
\par
Let $\A_N^{p,q}$  be the
sheaf of (complex valued) smooth $(p,q)$-forms on $N$,
$\dbar_N$ the tagential Cauchy-Riemann operator and $\Z_N^{p,q}$, (resp.
$\B_N^{p,q}$) the sheaf of closed (resp. exact) $(p,q)$-forms.
As usual we denote by $H^{p,q}(N)=\Z_N^{p,q}(N)/\B_N^{p,q}(N)$
the cohomology groups of the $\dbar_M$ complex on smooth forms.
The Poincar\'e lemma  is valid for
$\dbar_N$ (see \cite{Jurchescu:1994}), thus the complex:
\[
0\to\Omega^p_N\to\A_N^{p,0}\dbarNto\dots\dbarNto\A_N^{p,q}\dbarNto\dots
\]
is a fine resolution of $\Omega^p_N$. 
This implies that $H^{p,q}(N)\simeq H^q(\Omega_N^p)$.
\par
Let $\Eb_N$ be a homogeneous \CR vector bundle on $N$ (\ie a complex
 vector bundle with transition functions that are \CR) with fiber $E$,
and let $\E_N$ be the sheaf of its \CR sections.
\par
We denote by $\Eb_N^p$ the bundle
of \CR, $\Eb_N$-valued, $p$-forms,
with associated sheaf of \CR sections
$\E^p_N=\Omega^p_N\otimes_{\O_N}\E_N$.
\par
Let $\A^{p,q}_{N,\Eb_N}=\A_N^{p,q}\otimes_{\O_N}\E_N$, denote  by
$\dbar_{\Eb_N}$ the tangential Cauchy-Riemann operator on $\Eb_N$ and
let $\Z^{p,q}_{N,\Eb_N}$ (resp.\ $\B^{p,q}_{N,\Eb_N}$) be the sheaf of
$\dbar_{\Eb_N}$-closed (resp.\ exact) smooth forms with values in
$\Eb_N$. 
\par
Then $H^{p,q}(N,\Eb_N)=\Z^{p,q}_{N,\Eb_N}(N)/\B^{p,q}_{N,\Eb_N}(N)$,
but we also have:
\[
H^{p,q}(N,\Eb_N)=H^q(\E^p_N)=H^{0,q}(N,\Eb^p_N)=\Z^{0,q}_{N,\Eb^p_N}(N)/
\B^{0,q}_{N,
\Eb^p_N }(N).
\]
\par
For any open set $U\subset N$, the spaces $\A^{p,q}_{N,\Eb_N}(U)$ and
$\Z^{p,q}_{N,\Eb_N}(U)$ are Fr\'echet spaces with the topology of
uniform convergence of all derivatives on compact sets.
If $\B^{p,q}_{N,\Eb_N}(N)$ is closed in $\Z^{p,q}_{N,\Eb_N}(N)$, then
$H^{p,q}(N,\Eb_N)$, with the quotient topology, is also a Fr\'echet
space.

\section{Statements and proofs}
Let $\Eb_F$ be a $\L'$-homogeneous holomorphic vector bundle on $F$.
The $\L'$ action induces a natural $\L'$ action on
$\A_{F,\Eb_F}^{p,q}$, hence on
$H^{p,q}(F,\Eb_F)$, because the action of $\L'$ preserves closed and
exact forms.
Since $F$ is a compact complex manifold, $H^{p,q}(F,\Eb_F)$ is finite
dimensional and we can construct the $\K$-homogeneous holomorphic
vector bundle on $X'$:
\[
 \Hb_{X'}^{p,q}(F,\Eb_F)=\K\times_{\L'}H^{p,q}(F,\Eb_F)\text.
\]
\par
In a similar way we define a $\K_0$-homogeneous complex vector
bundle on $M'$: 
\[
\Hb_{M'}^{p,q}(F,\Eb_F)=\K_0\times_{\K_+'}H^{p,q}(F,\Eb_F)\text.
\]
\par
The following thorem has been proved by Le Potier (\cite{LePotier:1975},
see also \cite{FW:1984}):
\begin{thm}\label{thm:LePotier}
Let $X$, $X'$ and $F$ be as above, $\Eb_X$ a $\K$-homogeneous
holomorphic vector bundle on $X$ and $\Eb_X|^{}_F$ its restriction to
$F$.
Then there
exists a spectral sequence $^pE^{s,t}_r$, converging to
$H^{p,q}(X,\Eb_X)$, with
\[
^pE_2^{s,t}=\bigoplus_i
H^{i,s-i}\big(X',\Hb_{X'}^{p-i,t+i}(F,\Eb_X|^{}_F)\big).\qed
\]
\end{thm}
For $p=0$ the spectral
sequence collapses at $r=2$ and, recalling that $X'$ is a Stein
manifold, we obtain:
\[
H^{0,q}(X,\Eb_X)=H^{0,0}\big(X',\Hb_{X'}^{0,q}(F,\Eb_X|^{}_F)\big).
\]
Recalling that $H^{p,q}(X,\Eb_X)=H^{0,q}(X,\Eb_X^p)$ we finally obtain:
\begin{prop}\label{prop:complex}
Let $X$, $X'$ and $F$ be as above, $\Eb_X$ a $\K$-homogeneous
holomorphic
vector bundle on $X$ and $\Eb_X|^{}_F$ its restriction  to
$F$.
Then:
\[
H^{p,q}(X,\Eb_X)=H^{0,0}\big(X',\Hb_{X'}^{0,q}(F,\Eb^p_X|^{}_F)\big)
\]
as Fr\'echet spaces.\qed
\end{prop}
\par
A statement analogous to the last proposition  holds for $M$:
\begin{prop}\label{prop:cr}
Let $M$, $M'$ and $F$ be as above, $\Eb_M$ a $\K_0$-homogeneous \CR
vector bundle on $M$ and $\Eb_M|^{}_F$ its restriction  to
$F$.
Then:
\[
H^{p,q}(M,\Eb_M)=H^{0,0}\big(M',\Hb_{M'}^{0,q}(F,\Eb^p_M|^{}_F)\big)
\]
as Fr\'echet spaces.
\end{prop}
\begin{proof}
Fix $p$, $q$, let $\Z_{M'}=\pi_*(\Z_{M,\Eb^p_M}^{0,q})$,
$\B_{M'}=\pi_*(\B_{M,\Eb^p_M}^{0,q})$ and $\H_{M'}$ be the sheaf of
sections of $\Hb_{M'}^{0,q}(F,\Eb^p_M|^{}_F)$
We already know that
$H^{p,q}(M,\Eb_M)\simeq\Z_{M'}(M')/\B_{M'}(M')$.
\par
We now define a map
$\phi\colon\Z_{M'}\to\H_{M'}$ as follows.
\par
Let $U\subset M'$, $x\in U$, $x=g\K'_+$, $g\in\K_0$ and
$\xi\in\Z_{M'}(U)$.
Let $\xi_g=(g^{-1}\cdot \xi)|^{}_F$.
This is a closed $\Eb^p_M$-valued $(0,q)$-form on $F$, that determines a
cohomology class $[\xi_g]$  in $H^{0,q}(F,\Eb^p_M|^{}_F)$. 
Then  the class of $(g,[\xi_g])$ in
$\Hb_{M'}^{0,q}(F,\Eb^p_M|^{}_F)$ does not depend on the particular
choice of $g$, but only on $x$, hence it defines a section
$s_\xi=\phi(\xi)$ of $\H_{M'}$ on $U$.
\par
The sheaves $\Z_{M'}$, $\B_{M'}$, $\H_{M'}$, $\ker \phi$ are
$\O_{M'}$-modules and $\phi$ is a morphism of $\O_{M'}$-modules.
Since $\O_{M'}$ is fine, to prove that $\phi(M')$ is
continuous, surjective and
that its kernel is exactly $\B_{M'}(M')$ it sufices to check that
this is true locally, and this reduces to a straightforward
verification.
\end{proof}
\par
We prove now the main theorem of this paper:
\begin{thm}\label{thm:main}
Let $M$ and $X$ be as above, $\Eb_X$ a $\K$-homogeneous
holomorphic
vector bundle over $X$. Then:
\[
H^{p,q}(M,\Eb_X|^{}_M)\simeq\O_M(M)\otimes^{}_{\O_X(X)}H^{p,q}(X,\Eb_X)
\text.
\]
\end{thm}
\begin{proof}
Define $M'$ and $X'$ as above, 
fix integers $p,q\geq 0$ and let
$\Hb_{X'}=\Hb_{X'}^{0,q}(F,\Eb^p_X|^{}_F)$,
$\Hb_{M'}=\Hb_{M'}^{0,q}(F,\Eb^p_X|^{}_F)=\Hb_{X'}|^{}_{M'}$.
By Propositions \ref{prop:complex} and \ref{prop:cr}, we have that
$H^{p,q}(M,\Eb_X|^{}_M)=\Gamma(M',\Hb_{M'})$ and
$H^{p,q}(X,\Eb_X)=\Gamma(X',\Hb_{X'})$.
\par
Since $\dim_{\R}M'=\dim_{\C}(X')$, a global section of
$\Hb_{X'}$ that is zero on $M'$ must be zero on $X'$, \ie the
restriction
map $\Gamma(X',\Hb_{X'})\to\Gamma(M',\Hb_{M'})$ is injective.
\par
On the other hand, $X'$ is Stein, thus $\Hb_{X'}$ is generated at every
point by its global sections. Together with the fact that
$\Hb_{M'}=\Hb_{X'}|^{}_{M'}$ this implies that
\[
\Gamma(M',\Hb_{M'})=\O_{M'}(M')\otimes^{}_{\O_{X'}(X')}\Gamma(X',\Hb_{X'
} )
\text,
\]
where in the right hand side we implicitly identify global holomorphic
sections on $X'$ with
their restrictions to $M'$.
\par
The theorem follows from the observation that 
$\O_{M}(M)\simeq\O_{M'}(M')$ and $\O_{X}(X)\simeq\O_{X'}(X')$ because
the fiber $F$ of $\pi$ is a compact connected complex manifold.
\end{proof}
This, together with Corollary \ref{cor:dense}, implies the following.
\begin{cor}
With the same assumptions, the inclusion $j\colon M\to X$ induces a map:
\[
j^*\colon H^{p,q}(X,\Eb_X)\longrightarrow H^{p,q}(M,\Eb_X|^{}_M)
\]
that is continuous, injective and has
a dense range.\qed
\end{cor}

\section{An example}
Let $\G=\gr{SL}(4,\C)$, and $\Q$ be the parabolic subgroup of upper
triangular matrices.
We consider the real form $\G_0=\gr{SU}(1,3)$, identified with the
group of linear trasformations of $\C^4$ that leave invariant 
the Hermitian form associated to the matrix 
\[
B=\left(\begin{smallmatrix}0&0&0&1\\0&1&0&0\\0&0&1&0\\1&0&0&0
\end{smallmatrix}\right).
\]
Then $\G/\Q$ is the set of complete flags
$\{\ell^1\subset\ell^2\subset\ell^3\subset\C^4\}$ and $M=\G_0\cdot e\Q$
is the submanifold
$\{\ell^1\subset\ell^2\subset\ell^3\subset\C^4 \mid
\ell^3=(\ell^1)^\perp\}$.
\par
Let $\Q'$ be the set of block upper triangular matrices of the form
\[
\Q'=\left\{g=
\left(\begin{smallmatrix}*&*&*&*\\0&*&*&*\\
0&*&*&*\\0&0&0&*\end{smallmatrix}\right)
\mid g\in\G\right\},
\]
so that $M'$ is the totally real manifold
$\{\ell^1\subset\ell^3\subset\C^4\mid\ell^3=(\ell^1)^\perp\}$ and $M$
fibers over $M'$ with typical fiber $F$ isomorphic to $\C\P^1$.
The fibration is given by 
\[\tag{*}\label{ex1}
(\ell^1,\ell^2,\ell^3)\longmapsto(\ell^1,\ell^3).
\]
\par
Choose $\K$ to be the stabilizer in $\G$ of the subspaces
$V=\Span(\e_1-\e_4)$ and $W=\Span(\e_1+\e_4,\e_2,\e_3)$ so that $\K$ is
isomorphic to $\gr{S}(\gr{GL}(1,\C)\times\gr{GL}(3,\C))$ and $\K_0$ to
$\gr{S}(\gr{U}(1)\times\gr{U}(3))$.
Then $X$ is the set of flags
$\{\ell^1\subset\ell^2\subset\ell^3\subset\C^4\}$ in a generic position
with respect to the subspaces $V$ and $W$, and $X'$ is the set of flags
$\{\ell^1\subset\ell^3\subset\C^4\}$ in a generic position
with respect to $V$ and $W$. 
The map from $X$ to $X'$ given by \eqref{ex1} is a fibration with
typical fiber isomorphic to $\C\P^1$ and $X'$ is a Stein manifold.
\par
Finally let $\Eb=X\times\C$ be the trivial line bundle.
According to Propositions \ref{prop:complex} and \ref{prop:cr} the
cohomology of $M$ and $X$ is given by
\begin{align*}
H^{p,q}(X)&=H^{p,q}(X,\Eb)=H^{0,0}\big(X',\Hb_{X'}^{0,q}(F,\Eb^p|^{}_F)
\big),\\
H^{p,q}(M)&=H^{p,q}(M,\Eb|_M)
=H^{0,0}\big(M',\Hb_{M'}^{0,q}(F,\Eb^p|^{}_F ) \big).
\end{align*}
Recalling that $H^{p,q}(F)\simeq\C$ if $p=q=0$ or $p=q=1$ and $0$
otherwise, we obtain:
\[\begin{cases}
H^{p,q}(X)\simeq\O_X(X)\simeq\O_{X'}(X') &\text{if $p=q=0$ or
$p=q=1$;}\\
H^{p,q}(X)=0 &\text{otherwise;}
\end{cases}\]
and analogously:
\[\begin{cases}
H^{p,q}(M)\simeq\O_M(M)\simeq\O_{M'}(M') &\text{if $p=q=0$ or
$p=q=1$;}\\
H^{p,q}(M)=0 &\text{otherwise;}
\end{cases}\]
and it is clear that
\[
H^{p,q}(M)=\O_M(M)\otimes_{\O_X(X)}H^{p,q}(X).
\]

\def\MR#1{\href{http://www.ams.org/mathscinet-getitem?mr=#1}{MR\,#1}}
\providecommand{\bysame}{\leavevmode\hbox to3em{\hrulefill}\thinspace}

\end{document}